# An Improved Class Of Exponential Ratio- Type Estimator In The Presence Of Measurement Errors

Sachin Malik and †Rajesh Singh

Department of Statistics, Banaras Hindu University

Varanasi-221005, India

†Corresponding author

**Abstract**

In this paper we have suggested a family of estimators for the population mean in the presence of measurement errors. Expression for the mean squared error (MSE) of the suggested family is derived. An empirical study has been carried out to verify the theoretical findings.

## 1. Introduction

In survey sampling, the property of the estimators based on data usually presupposes that the observations are the correct measurements on characteristics being studied. Unfortunately, this ideal is not met in practice for a variety of reasons, such as non response errors, reporting errors, and computing errors. When the measurement errors are negligible small, the statistical inferences based on observed data continue to remain valid. On the contrary when they are not appreciably small and negligible, the inferences may not be simply invalid and inaccurate but may often lead to unexpected, undesirable and unfortunate consequences (see [1]). Some authors including [2], [2-9] have paid their attention towards the estimation of population mean $\mu_y$ of the study variable y using auxiliary information in the presence of measurement errors.

For a simple random sampling scheme, let $(x_i, y_i)$ be observed values instead of the true values $(X_i, Y_i)$ on two characteristics $(x, y)$ respectively for the $i^{th}$ (i=1.2….n) unit in the sample of size n.

Let the measurement errors be

$$u_i = y_i - Y_i \tag{1.1}$$

$$v_i = x_i - X_i \tag{1.2}$$

which are stochastic in nature with means zero and variances $\sigma_u^2$ and $\sigma_v^2$ respectively, and are independent. Further, let the population means of (x, y) be $(\mu_x, \mu_y)$, population variances of (x, y) be $(\sigma_x^2, \sigma_y^2)$ and $\sigma_{xy}$ and $\rho$ be the population covariance and the population correlation coefficient between x and y respectively (see [9]).

Let $w_u = \dfrac{1}{\sqrt{n}} \sum\limits_{i=1}^{n} u_i$, $w_y = \dfrac{1}{\sqrt{n}} \sum\limits_{i=1}^{n} (y_i - \mu_y)$

$w_v = \dfrac{1}{\sqrt{n}} \sum\limits_{i=1}^{n} v_i$, $w_x = \dfrac{1}{\sqrt{n}} \sum\limits_{i=1}^{n} (x_i - \mu_x)$

$C_x = \dfrac{\sigma_x}{\mu_x}$ and $C_y = \dfrac{\sigma_y}{\mu_y}$,

Also let

$k_1 = \bar{y} - \mu_y = \dfrac{1}{\sqrt{n}} (w_y - w_u)$

and, $k_2 = \bar{x} - \mu_x = \dfrac{1}{\sqrt{n}} (w_x - w_v)$, $\text{var}(\bar{y}) = \dfrac{\sigma_y^2}{n} \left[ 1 + \dfrac{\sigma_u^2}{\sigma_y^2} \right]$

and

$E(k_1) = E(k_2) = 0$

$E(k_1^2) = \dfrac{\sigma_y^2}{n} \left( 1 + \dfrac{\sigma_u^2}{\sigma_y^2} \right) = V_{ym}$

$E(k_2^2) = \dfrac{\sigma_x^2}{n} \left( 1 + \dfrac{\sigma_v^2}{\sigma_x^2} \right) = V_{xm}$

$$E(k_1 k_2) = \frac{\rho \sigma_y \sigma_x}{n} = V_{yxm}$$

In this paper, we have studied the behaviour of some estimators in presence of measurement error.

## 2. Some estimators in literature

The traditional ratio estimator is defined as

$$t_1 = \bar{y} \frac{\mu_x}{\bar{x}} \tag{2.1}$$

The bias and MSE of the estimator $t_1$ is given in [9], respectively, as

$$\text{Bias}(t_1) = \frac{1}{\mu_x} \left( R_m V_{xm} - V_{yxm} \right) \tag{2.2}$$

where $R_m = \frac{\mu_y}{\mu_x}$.

$$\text{MSE}(t_1) = \left[ \frac{\sigma_y^2}{n} + \frac{R_m^2 \sigma_x^2}{n} - \frac{2 R_m \rho_{yx} \sigma_y \sigma_x}{n} \right] + \frac{1}{n} \left[ \frac{\mu_y^2}{\mu_x^2} \sigma_v^2 + \sigma_u^2 \right] \tag{2.3}$$

MSE of the estimator $t_1$ can be re-written as

$$\text{MSE}(t_1) = M_{t_1}^* + M_{t_1} \tag{2.4}$$

where,

$M_{t_1}^* = \left[ \frac{\sigma_y^2}{n} + \frac{R_m^2 \sigma_x^2}{n} - \frac{2 R_m \rho_{yx} \sigma_y \sigma_x}{n} \right]$ is the mean squared error of the estimator $t_1$ without measurement error and

$M_{t_1} = \frac{1}{n} \left[ \frac{\mu_y^2}{\mu_x^2} \sigma_v^2 + \sigma_u^2 \right]$ is the contribution of measurement error in the estimator $t_1$.

[10] suggested an exponential ratio type estimator for estimating $\bar{Y}$ as

$$t_2 = \bar{y}\exp\left(\frac{\mu_x - \bar{x}}{\mu_x + \bar{x}}\right) \tag{2.5}$$

The bias and MSE of the estimator $t_2$ is given by [9], respectively, as

$$\text{Bias}(t_2) = \frac{1}{\mu_x}\left(\frac{3}{8}R_m V_{xm} - \frac{1}{2}V_{yxm}\right) \tag{2.6}$$

$$\text{MSE}(t_2) = \frac{\sigma_y^2}{n}\left[1 - \frac{C_x}{C_y}\left(\rho - \frac{C_x}{4C_y}\right)\right] + \frac{1}{n}\left[\frac{\mu_y^2}{4\mu_x^2}\sigma_v^2 + \sigma_u^2\right] \tag{2.7}$$

MSE of the estimator $t_2$ can be expressed as

$$\text{MSE}(t_2) = M_{t_2}^* + M_{t_2} \tag{2.8}$$

where,

$M_{t_2}^* = \dfrac{\sigma_y^2}{n}\left[1 - \dfrac{C_x}{C_y}\left(\rho - \dfrac{C_x}{4C_y}\right)\right]$ is the MSE of the estimator $t_2$ without measurement error, and

$M_{t_2} = \dfrac{1}{n}\left[\dfrac{\mu_y^2}{4\mu_x^2}\sigma_v^2 + \sigma_u^2\right]$ is the contribution of measurement error in the estimator $t_2$.

[11] suggested a regression type estimator $t_3$ as-

$$t_3 = \omega_1 \bar{y} + \omega_2(\mu_x - \bar{x}) \tag{2.9}$$

where $\omega_1$ and $\omega_2$ are constants.

[9] have derived the bias and MSE of the estimator $t_3$ respectively, given by

$$\text{Bias}(t_3) = \mu_y(\omega_1 - 1) \tag{2.10}$$

$$\text{MSE}(t_3) = \mu_y^2(\omega_1 - 1)^2 + \frac{1}{n}\omega_1^2\sigma_y^2 + \frac{1}{n}\omega_2^2\sigma_x^2 - \frac{2}{n}\omega_1\omega_2\rho\sigma_y\sigma_x + \frac{1}{n}\left(\omega_1^2\sigma_u^2 + \omega_2^2\sigma_v^2\right) \tag{2.11}$$

MSE of the estimator $t_3$ can be expressed as

$$\mathrm{MSE}(t_3) = M^*_{t_3} + M_{t_3} \tag{2.12}$$

where,

$$M^*_{t_3} = \mu_y^2(\omega_1 - 1)^2 + \frac{1}{n}\omega_1^2\sigma_y^2 + \frac{1}{n}\omega_2^2\sigma_x^2 - \frac{2}{n}\omega_1\omega_2\rho\sigma_y\sigma_x, \text{is the MSE of the estimator } t_2 \text{ without}$$

measurement error and

$$M_{t_3} = \frac{1}{n}\left(\omega_1^2\sigma_u^2 + \omega_2^2\sigma_v^2\right), \text{ is the contribution of measurement error in the estimator } t_2.$$

MSE of $t_3$ can be re-written as

$$\mathrm{MSE}(t_3) = (\omega_1 - 1)^2\mu_y^2 + \omega_1^2(a_1) + \omega_2^2 a_2 + 2\omega_1\omega_2(-a_3) \tag{2.13}$$

where,

$$a_1 = (V_{ym}), \qquad a_2 = (V_{xm}), \quad \text{and} \quad a_3 = (V_{yxm})$$

Optimising MSE of the estimator $t_3$ with respect to $\omega_1$ and $\omega_2$, we get the optimum values of $\omega_1$ and $\omega_2$ as

$$\omega_1^* = \frac{b_3 b_4}{b_1 b_3 - b_2^2} \quad \text{and} \quad \omega_2^* = -\frac{b_2 b_4}{b_1 b_3 - b_2^2} \tag{2.14}$$

where,

$$b_1 = \mu_y^2 + a_1, \qquad b_2 = -a_3, \qquad b_3 = a_2, \qquad \text{and} \quad b_4 = \mu_y^2.$$

Using these values of $\omega_1^*$ and $\omega_2^*$ from equation (2.14) into equation (2.13), we get the minimum MSE of the estimator $t_3$ as

$$\mathrm{MSE}(t_3)_{\min} = \left[\mu_y^2 - \frac{b_3 b_4^2}{b_1 b_3 - b_2^2}\right] \tag{2.15}$$

## [12] estimator under measurement error

[12] suggested an estimator $t_4$ as

$$t_4 = [m_1\bar{y} + m_2(\mu_x - \bar{x})]\exp\left(\frac{\mu_x - \bar{x}}{\mu_x + \bar{x}}\right) \tag{2.16}$$

Expanding equation (2.16) and subtracting $\mu_y$ from both sides, we get

$$(t_4 - \mu_y) = (m_1 - 1)\mu_y + m_1 k_1 - m_2 k_2 - \frac{k_2 m_1 R_m}{2} - \frac{k_1 k_2 m_1}{2\mu_x} + \frac{k_2^2 m_2}{2\mu_x} + \frac{3k_2^2 R_m}{8\mu_x} \tag{2.17}$$

On taking expectation of both sides of (2.18), we get the bias of the estimator $t_4$ to the order $O(n^{-1})$ as

$$\text{Bias}(t_4) = \frac{3R_m m_1 V_{xm}}{8\mu_x} - \frac{m_1 V_{yxm}}{2\mu_x} + \frac{m_2 V_{xm}}{2\mu_x} \tag{2.18}$$

From (2.17), we have

$$t_4 - \mu_y = \left[(m_1 - 1)\mu_y + m_1 k_1 - m_2 k_2 - \frac{k_2 m_1 R_m}{2}\right] \tag{2.19}$$

Squaring both sides of (2.19) and taking expectations, the MSE of the estimator $t_4$ to the order $O(n^{-1})$ is given by

$$\text{MSE}(t_4) = \mu_y^2 + m_1^2 A_1 + m_2^2 V_{xm} + m_1 m_2 A_2 - 2m_1 \mu_y^2 \tag{2.20}$$

where,

$$A_1 = \mu_y^2 + V_{ym} - R_m V_{yxm} + \frac{R_m^2 V_{xm}}{4} \quad \text{and} \quad A_2 = R_m V_{xm} - 2V_{yxm}$$

Mse($t_4$) will be minimum for

$$m_1 = \frac{2A_2 \mu_y^2}{A_2^2 - 4A_1 V_{xm}}$$

$$m_2 = \frac{4\mu_y^2 V_{xm}}{4A_1 V_{xm} - A_2^2}$$

(2.21)

## 4. The proposed estimator under measurement error

Following [10] and [12], we propose a modified class of estimator as

$$t_p = \left\{ m_{11} \frac{\bar{y}}{2} \left\{ \left(\frac{\bar{x}}{\mu_x}\right) \exp\left[\frac{(\bar{x} - \mu_x)}{(\bar{x} + \mu_x)}\right] + \left(\frac{\mu_x}{\bar{x}}\right) \exp\left[\frac{(\mu_x - \bar{x})}{(\mu_x + \bar{x})}\right] \right\} + m_{12}[\mu_x - \bar{x}] \right\} \exp\left[\frac{(\mu_x - \bar{x})}{(\mu_x + \bar{x})}\right] \quad (4.1)$$

Expanding equation (4.1) and then subtracting $\mu_y$ from both side

$$(t_p - \mu_y) = \left\{ (m_{11} - 1)\mu_y - m_{11}k_1 + \frac{9m_{11}R_m k_2^2}{8\mu_x} - m_{12}k_2 - \frac{m_{11}R_m k_2}{2\mu_x} \right.$$

$$\left. - \frac{m_{11}k_1 k_2}{2\mu_x} + \frac{m_{12}k_2^2}{2\mu_x} + \frac{3k_2^2 R_m}{8\mu_x} \right\} \quad (4.2)$$

On taking expectations of both sides of (4.2), we get the bias of the estimator $t_p$ to the order $O(n^{-1})$ as

$$\text{Bias}(t_p) = (m_{11} - 1)\mu_y + \frac{9m_{11}R_m V_{xm}}{8\mu_x} - \frac{m_{11}V_{yxm}}{2\mu_x} + \frac{m_{12}V_{xm}}{2\mu_x} + \frac{3V_{xm}R_m}{8\mu_x} \quad (4.3)$$

From equation (4.2), we have

$$(t_p - \mu_y) = (m_{11} - 1)\mu_y + m_{11}k_1 - m_{12}k_2 - \frac{m_{11}R_m k_2}{2\mu_x} \quad (4.4)$$

Squaring both sides of (4.2) and taking expectations, we get the MSE of the estimator $t_3$ to the order $O(n^{-1})$ as

$$\text{MSE}(t_p) = \mu_y^2 + m_{11}^2 q_1 + m_{12}^2 V_{xm} - 2m_{11}m_{12}V_{yxm} - 2m_{11}\mu_y^2 \quad (4.5)$$

where,

$$q_1 = \mu_y^2 + V_{ym} - R_m V_{yxm} + \frac{R_m^2 V_{xm}}{4}$$

Mse will be minimum for

$$\left. \begin{array}{l} m_{11} = \dfrac{\mu_y^2 V_{xm}}{q_1 V_{xm} - V_{yxm}^2} \\ m_{12} = \dfrac{\mu_y^2 V_{yxm}}{V_{yxm}^2 - q_1 V_{xm}} \end{array} \right\}$$

## 5. Empirical Study

**Data statistics:** The data used for empirical study has been taken from [13]

Where, $Y_i$ = True consumption expenditure,

$X_i$ = True income,

$y_i$ = Measured consumption expenditure,

$x_i$ = Measured income.

From the data given we get the following parameter values:

**Table 5.1: Data statistics**

| N | $\mu_y$ | $\mu_y$ | $\sigma_y^2$ | $\sigma_x^2$ | $\rho$ | $\sigma_u^2$ | $\sigma_v^2$ |
|---|---------|---------|--------------|--------------|--------|--------------|--------------|
| 10 | 127 | 170 | 1278 | 3300 | 0.964 | 36.00 | 36.00 |

**Table 5.2: Showing the MSE and PRE'S of the estimators with and without measurement errors**

| Estimators | MSE without meas. Error | Contribution of meas. Error in MSE | MSE with meas. Error | PRE'S |
|---|---|---|---|---|
| $\bar{y}$ | 127.800 | 3.600 | 131.400 | 100 |
| $t_1$ | 16.181 | 5.609 | 21.790 | 603.011 |
| $t_2$ | 25.925 | 4.102 | 30.027 | 437.596 |
| $t_3$ | 9.000 | 4.896 | 13.882 | 946.540 |
| $t_4$ | 8.972 | 4.932 | 13.904 | 944.943 |
| $t_p$ | 7.625 | 4.732 | 12.357 | 1063.365 |

**Conclusion**

From the Table 5.2, we conclude that the MSE of the proposed estimator $t_p$ is minimum as compared to the usual estimator $\bar{y}$ and the other considered estimators. We observe that the usual estimator $\bar{y}$ is least affected by the measurement errors.